\documentclass[12pt]{amsart}

\usepackage{amsmath,amssymb,enumerate}

\usepackage{epsfig,fancyhdr,color}

\usepackage{amssymb}
\usepackage{amsmath,amsthm}
\usepackage{bbm}
\usepackage{latexsym}
\usepackage{amscd}
\usepackage{psfrag}
\usepackage{graphicx}
\usepackage[all]{xy}

\usepackage[font=small,labelfont=bf]{caption}

\usepackage{mathtools}

\usepackage{tikz-cd}

\usepackage{stackrel}

\usepackage{mathabx,epsfig}

\usepackage{amsthm}

\usepackage[utf8]{inputenc}

\numberwithin{equation}{section}

\definecolor{NoteColor}{rgb}{1,0,0}

\newtheorem{theorem}{\rm\bf Theorem}[section]

\newtheorem{lemma}[theorem]{\rm\bf Lemma}
\newtheorem{corollary}[theorem]{\rm\bf Corollary}
\newtheorem*{theorem 1}{\rm\bf Proposition 1}
\newtheorem*{theorem 2}{\rm\bf Proposition 2}

\theoremstyle{definition}

\theoremstyle{remark}
\newtheorem{remark}[theorem]{\rm\bf Remark}

\newtheorem{example}[theorem]{\rm\bf Example}

\def\interieur#1{\mathord{\mathop{\kern 0pt #1}\limits^\circ}}

\begin{document}

\title{Music of moduli spaces}
\
\author{Robert Penner}
\address {\hskip -2.5ex Institut des Hautes \'Etudes Scientifiques\\
35 route des Chartres\\
Le Bois Marie\\
91440 Bures-sur-Yvette\\
France\\
{\rm and}~Mathematics Department,
UCLA\\
Los Angeles, CA 90095\\USA}
\email{rpenner{\char'100}ihes.fr}

 \date{\today}


\begin{abstract}
A musical instrument, the plastic hormonica, is defined here as a birthday present
for Dennis Sullivan, who pioneered and helped popularize the hyperbolic geometry underlying
its construction.
This plastic hormonica is based upon the Farey tesselation of the Poincar\'e disk
decorated by its standard osculating horocycles centered at the rationals.
In effect, one taps or holds points of another tesselation $\tau$ with the same decorating horocycles to produce sounds
depending on the fact that the lambda length of $e\in\tau$ with this decoration is always an integer.
Explicitly, tapping a decorated edge $e\in\tau$ with lambda length $\lambda$ produces
a tone of frequency $440\, \xi^{\lambda-12N}$, where $\xi^{12}=2$ and $N$ is some fixed positive integer shift of octave.
Another type of tap on edges of $\tau$ is employed to apply flips, which may be equivariant for a Fuchsian group preserving $\tau$.
Sounding the frequency for the edge after an equivariant flip, one can thereby audibly experience
paths in Riemann moduli spaces and listen to mapping classes.
The resulting chords, which arise from an ideal triangle complementary to $\tau$ by sounding the frequencies
of its frontier edges, correspond to a generalization of the classical Markoff
triples, which are precisely the chords that arise from the once-punctured torus.  In the other direction, one can
query the genera of specified musical pieces.
\end{abstract}

\maketitle

{\let\thefootnote\relax\footnote{{
{Keywords: Riemann moduli spaces, mapping class groups, lambda lengths, Penner coordinates, ideal triangulations, punctured surfaces, Farey tesselation, music.}
}}}
{\let\thefootnote\relax\footnote{{
{Note added in proof: Aaron Fenyes has developed the first implementation here {\tt https://vectornaut.github.io/harmonica} with the name changed to the more mellifluous and less
spell-check susceptible moniker {\sl horomonica}.}
}}}

\section*{Introduction}

Let{\let\thefootnote\relax\footnote{{
{It is a pleasure to thank Yi Huang for useful discussions and contributions, most especially in the proof of Theorem \ref{yay}, and
the referees for helpful comments.}
}}}
 ${\mathbb D}\subseteq{\mathbb C}$ denote the Poincar\'e disk lying in the complex plane
with its frontier {\it circle ${\mathbb S}^1\subseteq{\mathbb C}$ at infinity}.  The {\it modular group}  $${\rm PSL}_2={\rm PSL}_2({\mathbb Z}),$$ comprised of two-by-two integral unimodular matrices modulo multiplication by minus one,
acts isometrically on ${\mathbb D}$ preserving ${\mathbb S}^1$.  A {\it horocycle} is a curve in
${\mathbb D}$ of constant geodesic curvature unity, which is asymptotic to a unique point in
${\mathbb S}^1$ called its {\it center}, or equivalently, a Euclidean circle in ${\mathbb D}$
tangent to ${\mathbb S}^1$ at its center,  for example, the horocycle $h_*$ centered at $-1\in {\mathbb S}^1$ 
containing the origin $0\in{\mathbb D}$ of ${\mathbb C}$.

The {\it lambda length} $\lambda=\lambda(h,h')$ of a pair $h,h'$ of horocycles with distinct
centers is defined by  $$\lambda={\rm exp}~\delta/2,$$ where $\delta$ is the signed hyperbolic distance between 
$h$ and $h'$ along the geodesic connecting their centers, taken with a positive sign if and only if the horocycles are disjoint.     
A small calculation given in Section \ref{lam} shows that in fact, {\sl the lambda length of any two horocycles in the orbit ${\rm PSL}_2 (h_*)$ of $h_*$ is a positive integer}.  

Our basic idea here is to capitalize upon this integrality in order to define a virtual musical instrument, which is capable
of both human play and automatic actuation, in order to provide not only a mechanism for musical composition but also
a tool for the auditory probing of Riemann moduli spaces and mapping class groups of punctured surfaces.  
After all, human auditory acuity far exceeds the other senses,
so such a musical probe could help elucidate both geometry and algebra.

Recall that there are 12 musical tones in an octave, and frequency doubling corresponds to going up in tone by one octave.
The difference between two consecutive notes is called a {\it hemitone}.
The frequencies of musical notes in so-called equally tempered tuning, which is often implemented for instance
in electric pianos, are given by
$\omega_n=a\,\xi^n$, where $\xi=2^{1\over 12}=1.05946...$ and $a$ is a constant,
conventionally taken as $a=440$ Hz corresponding to the nearest note A above middle C, called ${\rm A}_4$, where the subscript
determines the octave.
The note with frequency $\omega_n$  is the musical note that is $n$ hemitones above or below ${\rm A}_4$.

Integrality of lambda lengths is exploited here, at first blush anyway, in the simplest way possible: the frequency assigned to a pair of horocycles with lambda length $\lambda$ is given by $440\,\xi^{\lambda-12N}=27.5\,\xi^\lambda$,
where the shift $N=4$ in octave is chosen so that small positive values of $\lambda$ lie just in the audible range, which for newborn humans is something like 20-20,000 Hz.  Thus, $N=4$ does the trick, with ${\rm A}_0$ of frequency 27.5 the effectively unhearable
note A four octaves below ${\rm A}_4$

More specifically, 
the {\it Farey tesselation} $\tau_*$ of ${\mathbb D}$, to be recalled in some detail the next section, is the orbit
$\tau_*={\rm PSL}_2 (e_*)$, where $e_*$ is the geodesic in ${\mathbb D}$
with endpoints $\pm 1\in{\mathbb S}^1$.  A collection of horocycles, one centered at each endpoint of a family of geodesics, is called a {\it decoration} on the family, and ${\rm PSL}_2 (h_*)$ thus provides a decoration
of $\tau_*$ called the {\it Farey decoration}.  Define the lambda length of a geodesic $e$ decorated by
respective horocycles $h,h'$ centered at its endpoints to be $$\lambda (e)=\lambda(h,h').$$

\begin{figure}[!h]
\begin{center}
\epsffile{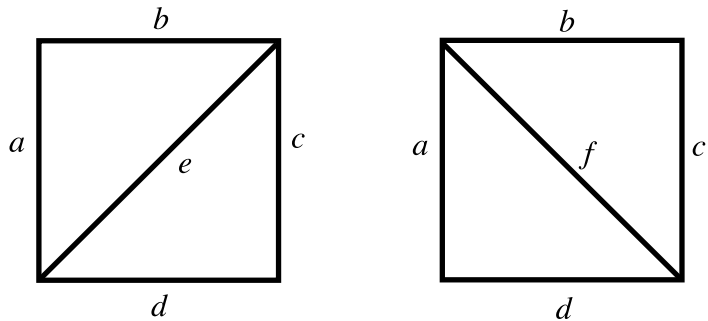}
\caption{Flip on an edge $e$.}
\label{fig:flip}
\end{center}
\end{figure}

A {\it flip} is the basic combinatorial move on a tesselation $\tau$, defined by removing a single edge $e$ from $\tau$, so as to produce a complementary ideal quadrilateral with diagonal $e$, and replacing $e$ in $\tau$ by the other diagonal $f$
of this quadrilateral, as depicted in Figure \ref{fig:flip}.  
The Ptolemy equation $$ef=ac+bd$$
describes the effect of flips on lambda lengths, which are here conflated with the edges themselves,
where $e$ and $f$ are diagonals of the quadrilateral with edges $a$ opposite $c$ and $b$ opposite $d$.

Thus, beginning with the {\sl untuned instrument} $\tau_*$, one might perform a sequence of flips
to {\sl tune} the instrument to another tesselation $\tau$, still with the original decoration, in order
to assign various integral lambda lengths, and hence this number of hemitones above some fixed tone, to the edges of $\tau$.

In a visual representation of $\tau$ on a computer screen not unlike a piano keyboard or a harp, one might thus touch the edges of $\tau$ in order to generate a musical sound at the frequency determined by the lambda length.  (We actually do something a little more elaborate than this to produce sounds in order to capture standard Western chords, cf. Section \ref{plastic1}.)

To probe moduli spaces, we consider a torsion-free subgroup $\Gamma$ of finite-index in ${\rm PSL}_2$,
so ${\mathbb D}/\Gamma$ is a {\it punctured arithmetic surface} 
$$F=F_g^s={\mathbb D}/\Gamma$$ of some genus $g\geq 0$ with a finite number $s\geq 1$
of punctures of negative Euler characteristic $2g-2+s<0$, and $\Gamma$ is a {\it punctured arithmetic surface group}.

Instead of flipping a single edge $e\in\tau$ in a $\Gamma$-invariant tesselation $\tau$, we perform flips 
one at a time on the edges in an entire orbit $\Gamma (e)$ in order to derive
another $\Gamma$-equivariant tesselation of ${\mathbb D}$, or in other words a triangulation of $F$ with its
vertices at the punctures, to be apprehended here as a $\Gamma$-equivariant retunings of $\tau$.  As explained in \cite{Pdec,Pbook}, sequences of these $\Gamma$-equivariant flips correspond to paths in the Riemann moduli space of $F$, and suitably periodic such sequences to elements of the mapping class group of $F$.  By listening to the sequences of corresponding tones, we might probe moduli spaces and mapping class groups alike.  

On the other hand, we might also take a fixed piece of music, and ask for the smallest genus of 
${\mathbb D}/\Gamma$ for which the piece may be played on a $\Gamma$-equivariant tuning of $\tau_*$.
As we shall see in Section \ref{plastic3}, any tune on one instrument using notes from only one octave, has genus at most three, for example,
the classic melody Happy Birthday to You, for Dennis Sullivan on this happy occasion of his 80th.

Dennis and I have been friends since the early 1980s.  He was lecturing in Boston on the measurable Riemann mapping theorem in the days before I married and received my doctorate, and we first met at a reception to celebrate the former.  During the next several years, we met often on Tuesdays, since I availed myself of his offer to all of us in the nearby Thurston gang to attend his Einstein seminar at CUNY and stay for dinner on his dime.  At times simultaneously raising families of roughly comparable ages, we were
at points very close indeed.  

Over the next many years, I would spend May with him at the IHES, and we would often meet elsewhere as well. In one of my Parisian visits, we had without any doubt the most efficient conversation of my life:  In the midst of doing math one evening after dinner, he abruptly looked me squarely in the eye and uttered ``math, love, children, wine and food," with each word lifting another finger on his hand, to which I responded ``yes, and music."  

Nevertheless, I offer here to Dennis for his birthday present a musical instrument, the {\sl plastic hormonica}, a precursor to which, called the {\sl hormonica}, was described in a footnote of \cite{Pbook}, which is reproduced here at the beginning of Section \ref{hormonica}.  I hope that Dennis among others will appreciate not only the math underlying the plastic hormonica, for instance the number-theoretic material on generalized Markoff triples, referred to here as {\sl triangular chords}, cf.~Section \ref{lam}, but also the tempting prospect of musically probing moduli spaces and mapping class groups. 

The oeuvre of Dennis pervades the constructions given here.
It has influenced my work, among so many others more generally in so many regards
and in so many directions, including his
leading one of the first seminars on hyperbolic geometry (at the IHES, in the 1980s).
He nurtured Bill Thurston, my own post-doctoral teacher, and gave all of us his
seminal works on the distribution of horocycles \cite{sullacta} in hyperbolic surfaces,
measures on the boundary of hyperbolic space \cite{sull2}, 
random walks in hyperbolic surfaces \cite{sull1}, the beautiful hyperbolic solenoid \cite{sull3}, 
among countless other contributions which have fundamentally impacted the development of  
hyperbolic space as a playground for all sorts of characters uncovering the richness of 
our modern hyperbolic geometry

 \section{Farey Tesselation}\label{farey}

Let us begin in the upper half-plane model 
$${\mathcal U}=\{x+iy\in{\mathbb C}:y>0\},$$ 
wherein horocycles are represented either as Euclidean circles tangent to ${\mathbb R}$ at their centers or as horizontal lines $\{y=c>0\}$
when the center is the point at infinity.
Let $h_n$ denote the horocycle with Euclidean diameter unity centered at $n\in{\mathbb Z}\subseteq{\mathbb R}$, for each $n\in{\mathbb Z}$ together with the horocycle $h_\infty=\{ y=1\}$.

\begin{figure}[!h]
\begin{center}
\epsffile{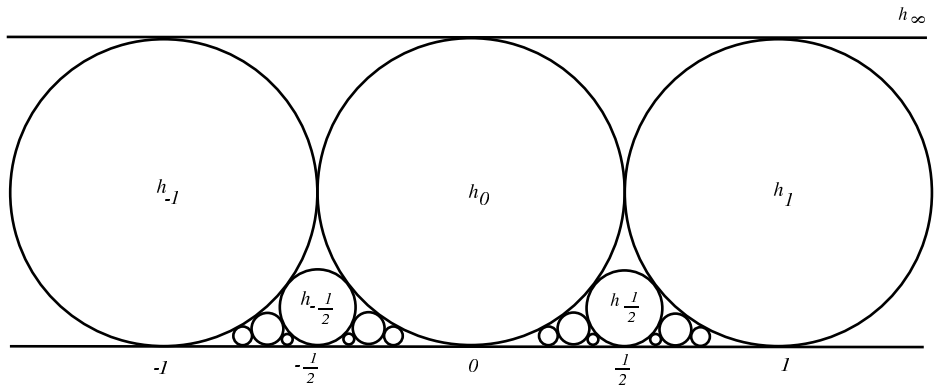}
\caption{The Farey decoration: a horocyclic packing of ${\mathcal U}$.}
\label{fig:6a}
\end{center}
\end{figure}

Two consecutive horocycles $h_n,h_{n+1}$ determine a triangular region bounded by the 
interval $[n,n+1]\subseteq{\mathbb R}$ together with the horocyclic segments connecting the horocycle centers to the point of tangency of $h_n$ and $h_{n+1}$.   There is a well-defined horocycle in each such triangular region which is tangent to $h_n$, to $h_{n+1}$, and to the real axis, and we let $h_{n+{1\over 2}}$ denote this horocycle, tangent to the real axis at the half-integer point $n+{1\over 2}$ and of Euclidean diameter ${1\over 4}$.
We may continue recursively  in this manner, adding horocycles tangent to the real axis and tangent to pairs of consecutive tangent horocycles, to produce a family of horocycles ${\mathcal H}$  in ${\mathcal U}$, cf. \cite{Ford}.  See Figure~\ref{fig:6a}.

\begin{figure}[!h]
\begin{center}
\epsffile{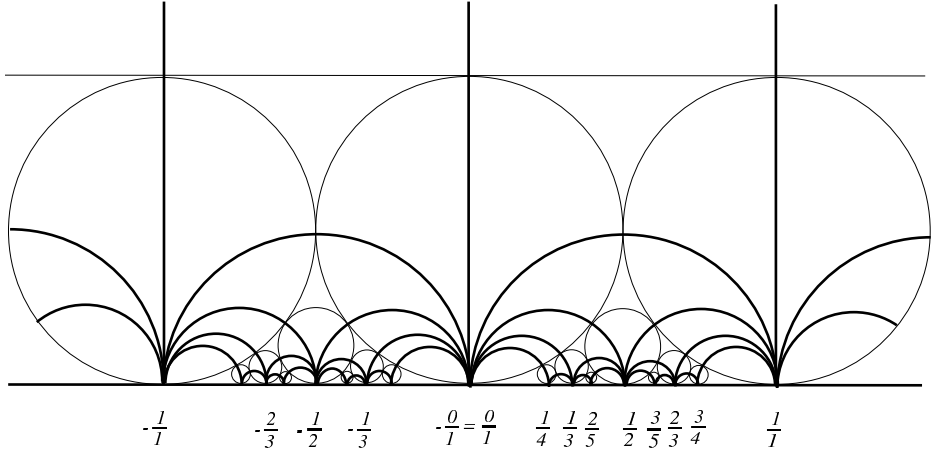}
\caption{Farey tesselation of ${\mathcal U}$.}
\label{fig:6b}
\end{center}
\end{figure}

\begin{lemma}\label{farey}
There is a unique horocycle in ${\mathcal H}$ centered at each extended rational point $\bar{\mathbb Q}={\mathbb Q}\cup\{\infty\}$, and the horocycle centered at ${p\over q}\in{\mathbb Q}$ has Euclidean diameter
${1\over{q^2}}$, where ${p\over q}$ is written in reduced form with $\infty =\pm{1\over 0}$.  Furthermore, the horocycles in ${\mathcal H}$ centered at distinct points ${p\over q},{r\over s}\in\bar{\mathbb Q}$ are tangent to one another if and only if $ps-qr=\pm 1$, and in this case, 
the horocycle in ${\mathcal H}$ tangent to these two horocycles is centered at ${{p+r}\over{q+s}}\in\bar{\mathbb Q}$.
\end{lemma}

It is easy to prove this lemma inductively starting with the second sentence.  
Now define the 
{\it Farey tesselation}  to be the collection
of hyperbolic geodesics in ${\mathcal U}$ that connect centers of tangent horocyles in ${\mathcal H}$; see Figure~\ref{fig:6b}.

\begin{figure}[!h]
\begin{center}
\epsffile{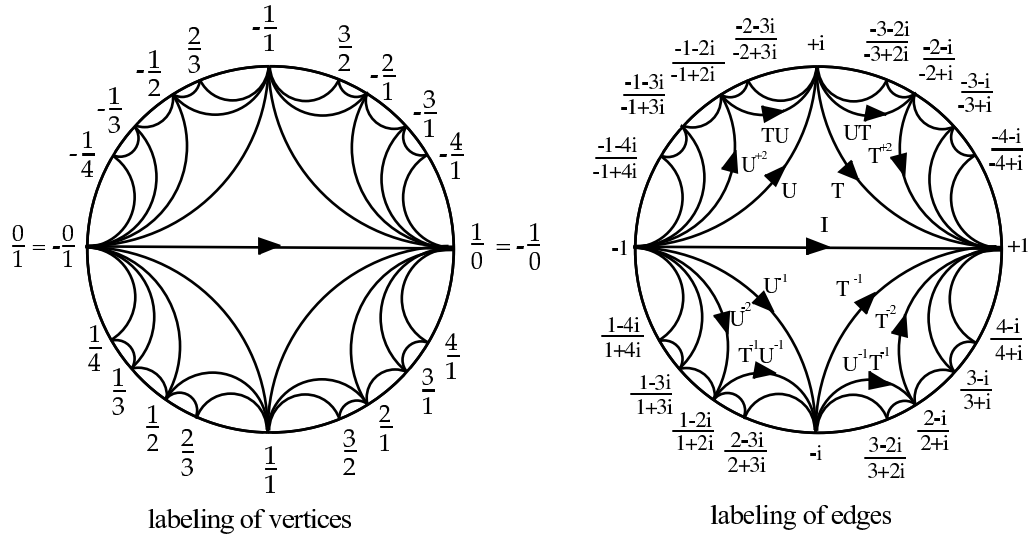}
\caption{Farey tesselation of ${\mathbb D}$.  Following Gauss,
we employ the right action on labels of edges.}
\label{fig:7}
\end{center}
\end{figure}

Finally define the {\it Farey tesselation} $\tau_*$ of ${\mathbb D}$ to be the image of the Farey tesselation of ${\mathcal U}$ under the Cayley transform $z\mapsto{{z-i}\over{z+i}}$
as illustrated in Figure \ref{fig:7}, where $\tau_*$ is regarded as a set of geodesics decomposing ${\mathbb D}$
into {\it ideal triangles}, i.e., regions bounded by three disjoint geodesics pairwise sharing ideal points at infinity. 
(The Cayley image of ${\mathcal H}$ in ${\mathbb D}$ is the orbit ${\rm PSL}_2(h_*)$ in the Introduction.)
The {\it generation} of a Farey point is the number of open ideal triangles met by the geodesic path from
the origin in ${\mathbb C}$ to the Farey point.

\begin{lemma}\label{modular}
The modular group (or rather its conjugate by the Cayley transform) leaves invariant the Farey tesselation $\tau _*$ and acts simply transitively on its oriented edges,
so oriented edges of $\tau_*$ are labeled by elements of ${\rm PSL}_2({\mathbb Z})$.
A generating set is given by any pair of
$$S=\left( \begin{matrix} 0&-1\\ 1&\hskip 1.5ex 0\end{matrix}\right ),
~~T=\left( \begin{matrix} 1&1\\ 0&1\end{matrix}\right ),
~~U=\left( \begin{matrix}  1&0\\ 1&1\end{matrix}\right ),$$
where $T^{-1}=SUS$ and $U^{-1}=STS$, and
a presentation in the generators $S,T$ is given by $S^2=1=(ST)^3$.
Reversal of orientation corresponds to precomposition with $S$.
\end{lemma}

In the sequel, we shall intentionally conflate a matrix 
$\begin{psmallmatrix}a&b\\c&d\end{psmallmatrix}$ with its
projectivization as well as with its corresponding linear fractional
transformation $z\mapsto {{az+b}\over{cz+d}}$ acting on $\mathcal U$,
and even its Cayley conjugate acting on ${\mathbb D}$.

\section{Lambda Lengths and Triangular Chords}\label{lam}

Direct calculations in the upper half-plane given in \cite{Pdec,Pbook} prove

\begin{lemma}\label{lemma1}
If $h,\bar h$ are two horocycles in ${\mathcal U}$ with respective centers $x,\bar x\in{\mathbb R}$ and Euclidean diameters
$\delta,\bar\delta>0$, then 
$$\lambda(h,\bar h)={{|x-\bar x|}\over{\sqrt{\delta\bar \delta}}},$$
while if $\bar h$ is centered at infinity with Euclidean height $H$, then we have $\lambda(h,\bar h)=\sqrt{{H\over{\delta}}}$. 
Furthermore, if $\gamma=\begin{psmallmatrix}a&b\\c&d\\\end{psmallmatrix}\in{\rm PSL}_2({\mathbb R})$ maps $h$ to $\bar h$, then 
$\bar x=\gamma(x)$ and 
$$\bar \delta =\gamma'(x)\delta={\delta\over{(cx+d)^2}},$$
while if $\gamma(x)=\infty$, then $\bar h$ has height ${\delta^{-1}\over{c^2}}$.
\end{lemma}

\begin{corollary} \label{cor:int} The lambda length $\lambda$ of the pair of horocycles in the Farey decoration 
with centers ${p\over q},{r\over s}$, where $p,q$ and $r,s$ are each coprime pairs, is given by
$\lambda=|ps-qr|$.
\end{corollary}

\begin{proof} For the horocycles of respective Euclidean diameters ${1\over q^2}$ and
${1\over s^2}$ centered at ${p\over q}<{r\over s}$, we find the lambda length
$${{|{p\over q}-{r\over s}|}\over\sqrt{{1\over{q^2}}{1\over{s^2}}}}=|ps-qr|=qr-ps ~{>0}.$$
\end{proof}

\noindent This proves the basic integrality property of lambda lengths on the Farey decoration which is at the heart of this paper.

Choose any three distinct rational numbers
${p\over q}<{r\over s}<{u\over v}$, each with coprime numerator and denominator, so the Farey decoration
assigns respective Euclidean diameters ${1\over {q^2}}$, ${1\over {s^2}}$, ${1\over{v^2}}$.
The triple of pairwise lambda lengths of these three horocycles is given by 
$$\{qr-ps, us-rv,uq-pv\},$$ and we call such a triple
a {\it (triangular) chord}.
Note that lambda lengths are invariant under the action of ${\rm PSL}_2({\mathbb R})$, so in particular the triple of lambda lengths is invariant under the diagonal action of the modular group on $\bar{\mathbb Q}^3$.

\begin{figure}[!h]
\begin{center}
\epsffile{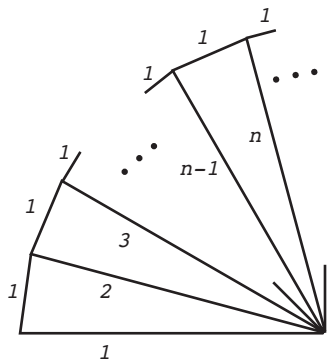}
\caption{The hyperfan $P$ with lambda lengths indicated.}
\label{fig:fan}
\end{center}
\end{figure}

\begin{example}\label{piano} For any $n\geq 1$, each of $\{ 1,n,n+1\} $, $\{ 1,n,2n+1\}$ and $\{ 1,n+1,2n+1\}$ is a triangular chord.
To see this in the notation of Figure \ref{fig:7}, serially flip along each of the edges $U, U^2,U^3, \ldots ,U^n,\ldots$
in this order to produce the polygon $P$ depicted in Figure \ref{fig:fan}, decomposed into a collection of triangles sharing a vertex with respective triples of lambda lengths $\{ 1,1,2\}, \{1,2,3\}, \{ 1,3,4\} ,\ldots , \{1,n-1,n\}, \dots$, proving the first part. (The weighted edges interior to $P$ comprise a ``hyperfan'' in the sense
of \cite{Puniv,Pbook}.)  Now flip on
edges in the frontier $P$ for the other two classes of chords. 
\end{example}

 Triangular chords generalize
the classical diophantine Markoff triples in the following sense.
Start from a triple of lambda lengths satisfying the Markoff equation 
$$x^2+a^2+b^2=3abx,$$ such as $x=a=b=1$.
The sum of the two roots of this quadratic equation in $x$ is given by $3ab$ according to the quadratic formula, and yet immediately $3ab=x+(a^2+b^2)/x$ for a solution to this equation.  It follows that if 
$(a,b,x)$ is a Markoff triple, i.e., a solution to the Markoff equation, then so too is
$(a,b,(a^2+b^2)/x)$, and we recognize this as the Ptolemy transformation applied
to the triple $(a,b,x)$ in the surface $F_1^1$.
In fact by well-known results \cite{Cassels},
the Markoff triples are exactly the triangular chords arising from any finite sequence of
$\Gamma$-equivariant flips on $\tau_*$, for any marked $\Gamma<{\rm PSL}_2$
with ${\mathbb D}/\Gamma$ homeomorphic to $F_1^1$.

First of all, one can relax this to consider $\Gamma<{\rm PSL}_2$ with any fixed topological type
and take finite sequences of $\Gamma$-equivariant flips to define various modular classes of triangular chords.  

\begin{remark} \label{jhc} A fundamental fact \cite{Pdec,Pbook}, essentially going back to J.H.C.~Whitehead, is that finite sequences of $\Gamma$-equivariant flips act transitively on the collection of all tesselations of ${\mathbb D}$ covering ideal triangulations of the topological surface ${\mathbb D}/\Gamma$.  It follows that the collection of $\Gamma$-equivariant tesselations of ${\mathbb D}$ arising from finite sequences of $\Gamma$-equivariant flips on any one of them depends only upon the topological type of the surface ${\mathbb D}/\Gamma$.
\end{remark}

Or one can ``go all the way" with the ``universal" definition above for triangular chords arising from any finite sequence of 
simple (i.e., non-equivariant) flips.  (See \cite{FP,Puniv,Pbook} for the sense in which the latter is universal.)

Let $gcd(n_1,\ldots , n_m)$ denote the greatest common divisor of any finite set $\{ n_1,\ldots ,n_m\}\subseteq {\mathbb Z}$ 
of integers.  The rest of this section is dedicated to the proof of the following theorem.

\begin{theorem} \label{yay} A triple $\lambda_i\in{\mathbb Z}$, for $i=1,2,3$ occurs as a triangular chord if and only if the following two
conditions hold whenever $\{i,j,k\}=\{ 1,2,3\}${\rm :}

\medskip

\leftskip .5cm

\noindent {\bf 1.} $gcd(\lambda_i,\lambda_j)$ divides $\lambda_k$;

\medskip

\noindent {\bf 2.} if $n=gcd( \lambda_1,\lambda_2,\lambda_3)$  is even, then some ${\lambda_i}/n$ is even.

\leftskip=0ex
\end{theorem}

\begin{lemma} \label{lemmaq}Given $\lambda_i\in{\mathbb Z}$, for $i=1,2,3$ satisfying Condition 1, let 
$n=gcd(\lambda_1,\lambda_2,\lambda_3)$.  Then
$n=gcd( \lambda_i,\lambda_j)$ as well, independently of the pair $i,j\in\{ 1,2,3\}$.
\end{lemma}

\noindent This means that the $\lambda_i'={\lambda_i/ n}$, for $i=1,2,3$ are pairwise coprime
(and hence trivially satisfy both conditions 1 and 2).  In particular, at most one of $\lambda_i'$, for $i=1,2,3$, can be even.

\begin{proof} [Proof of Lemma \ref{lemmaq}] Since $gcd( \lambda_i,\lambda _j)=n\, gcd( \lambda'_i,\lambda' _j)$
divides $\lambda_k$, for $\{ i,j,k\}=\{ 1,2,3\}$, it follows that $gcd(\lambda'_i,\lambda' _j)$ divides ${\lambda_k/ n}=\lambda_k'$.
Thus, $gcd( \lambda'_i,\lambda_j')$ divides each $\lambda_i$, for $i=1,2,3$, and hence also $n$.  Conversely, 
$gcd(\lambda_1',\lambda_2',\lambda_3')$ trivially divides $gcd(\lambda_i',\lambda_j')$, and hence finally, they must be equal.
\end{proof}

\begin{proof}[Proof of Theorem~\ref{yay}]
We start with necessity of Condition 1 and suppose that the triangular chord in question arises from the triple
${p\over q}<{r\over s}<{u\over v}$.  Choose $\gamma\in{\rm PSL}_2({\mathbb R})$ with $\gamma ({r\over s})={0\over 1}$
as well as $\gamma'({p\over q})=1=\gamma'({u\over v})$, and let 
$\gamma({p\over q})={{p'}\over {q'}}$, $\gamma({u\over v})={{u'}\over {v'}}$.
By the second part of Lemma \ref{lemma1}, the resulting horocycles centered at ${{p'}\over{q'}}$ and ${{u'}\over{v'}}$ have unchanged
diameters, and now the diameter of the horocycle at ${0\over 1}$ is given by some $\delta>0$.

The resulting lambda length of the horocycles
centered at ${p'  \over q'}$ and ${{u'}\over {v'}}$ is given by
$u'q'-p'v'$, which agrees with $uq-pv$ owing to ${\rm PSL}_2({\mathbb R})$-invariance of lambda lengths, while the other two 
lambda lengths are expressed as $-{{p'}\over\sqrt{\delta}}=qr-ps$ and ${{u'}\over\sqrt{\delta}}=us-rv$ by the first part of Lemma \ref{lemma1}, again using invariance of lambda lengths.

It follows that $uq-pv=u'q'-p'v'$ is a multiple of the greatest common divisor of $-p'$ and $u'$ and hence of
$-{{p'}\over\sqrt{\delta}}=qr-ps$ and ${{u'}\over\sqrt{\delta}}=us-rv$ as well.
Thus, Condition 1 is indeed necessary.
It follows that $\{ 10,12,15\} $ is not a triangular chord, for example.  

Continuing with the proof of Theorem \ref{yay}, let us henceforth adopt the abiding notation
$$A=\lambda _1, ~B=\lambda_2, ~C=\lambda_3, ~~{\rm and}~a=A/n, ~b=B/n, c=C/n,$$
where $n=gcd(A,B,C)$.

We turn to necessity of Condition 2 and suppose not, so that $n$ is even, and each of
$a,b,c$ is odd.  Transitivity of ${\rm PSL}_2$ on triples in $\overline{\mathbb Q}={\mathbb Q}\cup\{\infty\}$
allows us to position respective
vertices for the ideal triangle realizing  the triangular chord
$A,B,C$ at ${0\over 1}, A/r=na/r, B/s=nb/s$, for some $r,s\in{\mathbb Z}$.
It follows that
$$nc=C=As-Br=n(as-br),$$
where we have absorbed the absolute value into the signs of $r,s$.

Thus, $c=as-br$, where each of $a,b,c$ is assumed to be odd, so
$\{r,s\}$  must be equal to $\{ 0,1\}$ modulo 2.  This is a contradiction since
$na/r$ and $nb/s$ are each reduced fractions, with $n$ even and precisely
one of $r,s$ even, and this establishes necessity of Condition 2.

Now we suppose that $A,B,C$ is a triple satisfying both Conditions 1 and 2 and must
establish that $A,B,C$ is a chord.  
Write $$n=2^{e_0}~p_1^{e_1}\cdots  p_m^{e_m}$$
as a product distinct primes, where $e_0\geq 0$ and $e_1,\ldots ,e_m\geq 1$.

If $e_0= 0$, then it can be ignored, and the discussion below for the other prime factors pertains; otherwise, by Condition 1 and the fact already discussed that
$a,b,c$, are pairwise coprime, at most one
of these three can be even, say it is $c$.

Let us construct the triangular chord with vertices ${0\over 1}$, $a/r$, $b/s$ realizing $a,b,c$  where
$c=as-br$.  From the parities of the $a,b,c$, it follows
that $r$ and $s$ must have the same parities.

If both of $r,s$ are even, then $c$ may be rewritten as $$c=a(b+s)-b(a+r)$$ with
${0\over 1}$, $a/(a+r)$ and $b/(b+s)$ still in reduced form and yielding lambda lengths
$a,b,c$.  Replacing 
$$s\mapsto b+s~~{\rm and}~~r\mapsto a+r,$$ we may assume without loss
that $r$ and $s$ are each odd, and in this case, each of
${0\over 1}$, $(2^{e_0})a/r$, $(2^{e_0})b/s$ is a reduced fraction, thus
providing an ideal triangle for the chord $(2^{e_0})a, (2^{e_0})b, (2^{e_0})c$.

This shows that if $a,b,c$ is a chord, then so too must be its homothetic scaling
by $2^{e_0}$.  Now we proceed to include each of the remaining prime powers in $n$
by adjusting the denominators $r,s$ so that they become coprime to the numerator
we hope to introduce. Namely, we proceed by induction under the assumption that $r$ and $s$ are coprime to 
$$F=2^{e_0}p_1^{e_1}\cdots p_k^{e_k}$$ and modify $r,s$ so as to be coprime to
$2^{e_0}p_1^{e_1}\cdots p_{k+1}^{e_{k+1}}$.  It follows that 
$${0\over 1},~~(2^{e_0}\,  p_1^{e_1}\cdots p_{k+1}^{e_{k+1}})a/r,~~(2^{e_0}\, p_1^{e_1}\cdots p_{k+1}^{e_{k+1}})b/s$$ are already
in reduced form and achieve the desired lambda lengths.

Suppose that $p_{k+1}$ divides one of $r$ and $s$, and suppose without loss that it is $r$.
Note first that if $p_{k+1}$ were to also divide $s$, then it would moreover divide $(F)c=(F)as-(F)br$,
and hence $p_{k+1}$ would divide $c$.  By relative primality then $p_{k+1}$ cannot divide 
$(F)a$ or $(F)b$, and so $(F)b+s$ and $(F)a+r$ are each likewise 
coprime to $p_{k+1}$ and to $F$.

It follows that
$${0\over 1},~~(F)a/((F)a+r),~~(F)b/((F)b+s)$$
is a triple of reduced fractions realizing the triple of lambda lengths
$(F)a, (F)b, (F)c$.  Now take
$$s\mapsto (F)b+s~~{\rm and}~~r\mapsto (F)a+r$$ to
return to the case that $p_{k+1}$ and $F$ are coprime to both $r$ and $s$,
which was treated before.

It remains only to consider the case that $p_{k+1}$ divides $r$ but not $s$.  Since $p_{k+1}\neq 2$,
$p_{k+1}$ is coprime to both 
$(F)a+r$ and $2(F)a+r$ and likewise coprime to at least one of $(F)b+s$ and $2(F)b+s$. 

If $p_{k+1}$ is coprime to $(F)b+s$, then take 
$$s\mapsto (F)b+s~~{\rm and}~~r\mapsto (F)a+r,$$
and if $p_{k+1}$ divides  $(F)b+s$, then take 
$$s\mapsto 2(F)b+s~~{\rm and}~~r\mapsto 2(F)a+r .$$

In any case, the resulting $r$ and $s$ still satisfy $c=as-br$ and are respectively
coprime to $(F)a$ and $(F)b$ by the Euclidean Algorithm.  Furthermore, both are
coprime to $p_{k+1}$ and $F$, again reducing to the earlier case.

This finally provides an ideal triangle with vertices
$${0\over 1}, ~~(F p_{k+1}^{e_{k+1}})a/r=(2^{e_0}p_1^{e_1}\cdots p_{k+1}^{e_{k+1}})a/r,~~(F  p_{k+1}^{e_{k+1}})b/s=(2^{e_0}p_1^{e_1}\cdots p_{k+1}^{e_{k+1}})b/s$$
realizing lambda lengths
$$(2^{e_0}  p_1^{e_1}\cdots p_{k+1}^{e_{k+1}})a, (2^{e_0}  p_1^{e_1}\cdots p_{k+1}^{e_{k+1}})b, (2^{e_0} p_1^{e_1}\cdots p_{k+1}^{e_{k+1}})c.$$
The algorithmic procedure continues in this way through the remaining prime
factors of $n$.
\end{proof}

It follows from the theorem that pairwise coprime triples are triangular chords, as was noted before, and moreover the union of all triangular chords, as determined in the theorem,  is closed under integral homothety.

Furthermore, the theorem provides a generalization of the Bezout identity: For any pair $a,b$ of coprime integers and any non-zero integer $n$, there exists a pair of integers $r,s$ so that $as-br=1$ and $gcd(n,r)=gcd(n,s)=1$, namely, in the notation of the proof
either exactly one of $a,b,c$ is even or all three are odd.

\medskip

\leftskip=.3in\rightskip .3in

\noindent 

\bigskip

Basic problems and questions include:

\leftskip .5in\rightskip .3in

\noindent $\bullet$ In light of Remark \ref{jhc}, characterize the triangular chords for various topological types of punctured arithmetic surfaces.    
In particular, how does adding a single puncture $F_g^s\mapsto F_g^{s+1}$ affect the collection  of equivariant triangular chords?  

\leftskip=0ex\rightskip=0ex

\section{The hormonica and tempering}\label{hormonica}

Here is the footnote from \cite{Pbook} which defines the {\it hormonica}, the starting point for this paper.

\bigskip

\leftskip .2in\rightskip .2in

{\footnotesize At the risk of proving beyond any doubt that we have too much spare time, let us remark that there is a musical instrument, the ``hormonica'', \index{hormonica} based on this observation [{\it integrality of lambda lengths for the Farey decoration}] as follows.  Begin with the Farey tesselation $\tau_*$ as the untuned instrument regarded as
drawn before you on the computer screen.  Perform a sequence of flips by serially selecting edges to produce another tesselation $\tau$ of the Poincar\'e disk likewise displayed on the computer screen.  Choose some basic frequency, say middle C, to represent unity, so that any natural number may be interpreted as a multiple of this frequency.  In this way, each edge of  $\tau$ with its integral lambda length can be ``plucked'' to produce a corresponding tone.  Moreover, each triangle complementary to $\tau$ can be
``tapped'' to produce a triple of tones or chord.  Attributes such as 
duration or timbre could be introduced as further aspects of the tuning process.
Maybe this is crazy, but it could be fun.  On the other hand, it is difficult to probe the combinatorics of moduli spaces visually, and the analogous hormonicas  based upon
tesselations of a fixed surface could provide an auditory tool towards this end.}

\bigskip

\leftskip=0ex\rightskip=0ex

There are several difficulties with this definition.  The first problem is that frequencies of
musical notes are exponential (in equal-tempered tuning
given by $440\, \xi^n$, for integral $n$ with $\xi^{12}=2$) 
and not multiplicative in this sense and are rationally 
rather than integrally related, as discussed below.
Another more amusing limitation is that traditional musical chords, as given in Table \ref{table1},
cannot be triangular chords on the original hormonica,
since neither minor chords 10:12:15  nor diminished chords 160:192:231 (nor its sometimes useful approximation 20:24:29) satisfy Condition 1 of 
Theorem \ref{yay}.  

However, not all major, minor or diminished chords can arise from ideal triangles on the plastic hormonica either, with its system $\lambda\mapsto 27.5\, \xi^\lambda$ of assigning frequencies to lambda lengths; for example, the chords ${\rm B}_0$ major (2,6,9), ${\rm C}_1$ minor (3,6,10), and ${\rm B}_0$ diminished (2,5,8) violate Condition 1 of 
Theorem \ref{yay}, where the triples are numbers of hemitones of individual notes above ${\rm A}_0$, cf.~the first two columns of Table \ref{table1}.  This deficiency in the musical instrument will be rectified by activating tones with several fingers.
 
\begin{table}[tbhp]
\centering
\title{\bf Frequency Ratios of 3-Tone Musical Chords}

\begin{tabular}{ccccc}\\
$\underline{\rm Chord}$	&	$\underline{\rm Hemitones~between~notes}$		&		$\underline{\rm Ideal~Frequency~Ratios}$\\
Major	&	4-3	&		4:5:6	\\
Minor	&	3-4	&		10:12:15\\
Diminished	&	3-3	&		160:192:231\\
\hline
\end{tabular}
\caption{The tones comprising major, minor and diminished chords.}
\label{table1}
\end{table}


This table illustrates both the utilities and deficiencies of equally tempered tuning,
for example, the ideal major chord frequency ratios in equally tempered tuning
are not rational, but are approximated by $\xi^4=1.2599... \approx 5/4$ and $\xi^7=1.498... \approx 6/4$. 
But in what sense do these irrational frequency ratios fail to be ``ideal"?

\begin{figure}[!h]
\begin{center}
{\epsfysize3.8in\epsffile{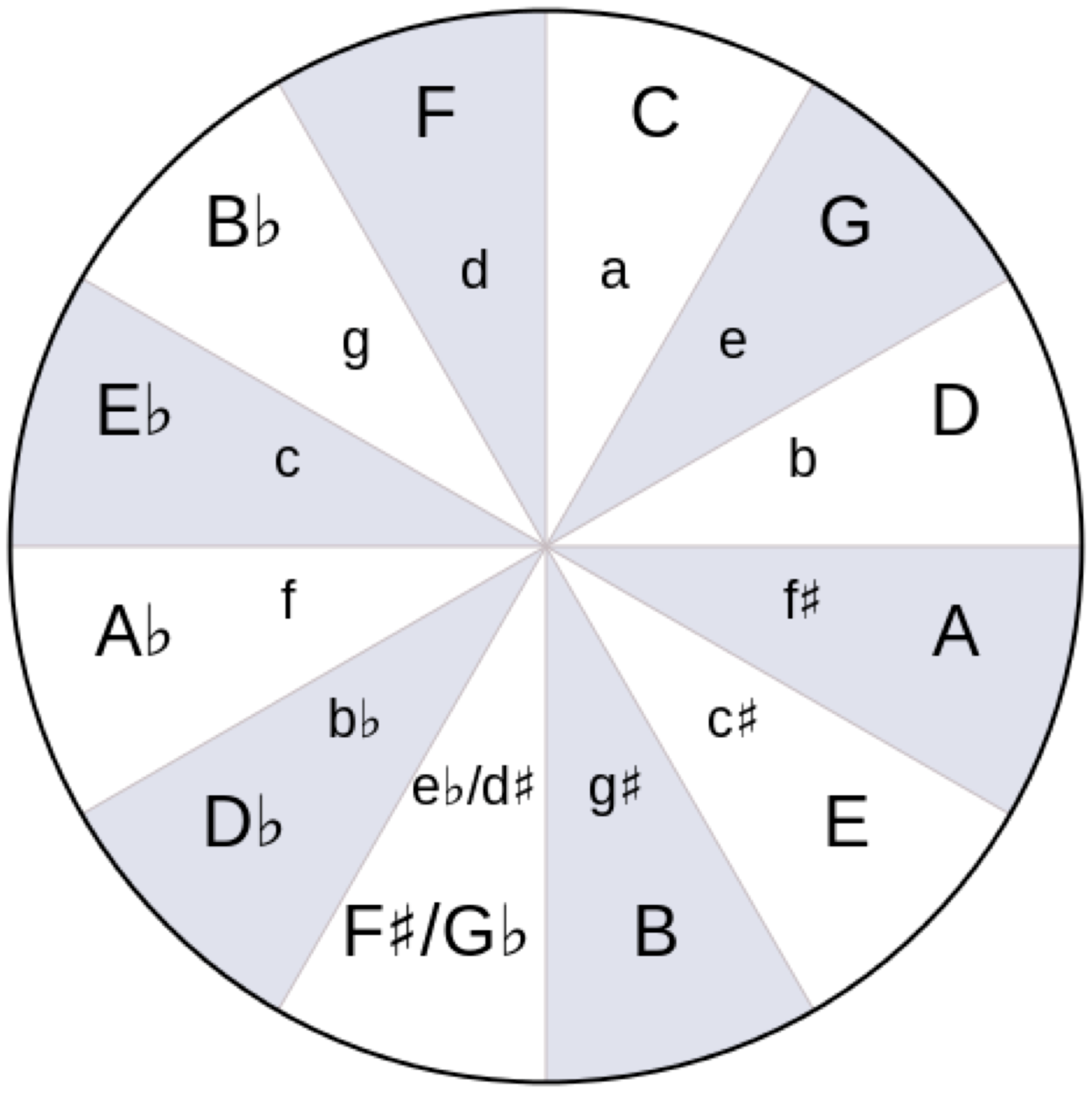}}
\vskip -.5in
\caption{The ideal musical notes arise from taking perfect fifths, namely frequencies in the ratio 3:2, starting from C and going clockwise around the outside of the so-called
circle of fifths illustrated here, so for instance C to G , G to D, and so on are all perfect fifths.  
Since $(3/2)^{12}=129.74633...$ is not a power of 2, there is a resulting inconsistency in octaves.  The lower case inner circle likewise give perfect fifths beginning from a.  Figure courtesy of Mysid/Andeggs.}
\label{fig:circle}
\end{center}
\end{figure}

That equal tempering is not the most harmonious tuning system goes back to the ancient Greeks, who argued that we humans find sonorous those pairs of frequencies whose ratio is a rational number with small denominator and numerator,
for instance the perfect fifth 3:2 and major third 6:5 in the ideal major chord; see Figure \ref{fig:circle}.  

A ``tempering'' is
a compromise in assigning frequencies to notes that respects the octave as well as preserves perfect fifths, on top of other desired rational frequency ratios within an octave, as least as much as possible according to certain criteria.
It is thus a compromise between the physics of a vibrating string
and the physiological perceptions of multiple frequencies.

Pythagoras among others has proposed such compromise systems of rationally related frequencies in an octave\footnote
{Both interesting and beautiful, here are the Pythagorean frequency ratios in the key of C: C(1:1), D(9:8), E(81:64), F(4:3), G(3:2), A(27:16), B(243:128), C(2:1).}.
An even more elaborate scheme is used to tune an ordinary piano across its 87 hemitones=7${1\over 3}$ octaves.
Moreover, triples of tones are harmonious when their integral frequencies share a 
common multiple in each of their eight smallest multiples according to Helmholtz.

The point is that books have been written \cite{tempering1,tempering0} on this topic of tempering with its rich history,
and it is more than we can discuss here.  We needed at least to give sense to the term, this
compromise in tuning away from equal tempering to achieve some specified collection of diophantine frequency ratios.

In practice for a stringed musical instrument, one ``fudges'' a little here and there to achieve
desired diophantine ratios at the expense of others, for instance in tuning a piano or in placing the frets of a guitar.
In a viola or violin, direct continuous tempering is determined in real time by the musician based on where the bow and finger meet
the string among other things, or in a guitar, where notes are determined not only by the fret but also by bending the string
to change its tension.
The plastic hormonica allows any type of specified tempering analogous to both a piano
and a guitar, as well as its real-time refinement as for a violin, at heart exponential based on equal temperament

\section{The plastic hormonica}\label{plastic}

This section is finally dedicated to the definition of a musical instrument, the {\it plastic hormonica}, and has three parts: first, 
the definition of the untuned instrument, which may provide an interesting interaction with hyperbolic geometry; second, the methods of tuning the instrument, which can be done either one ideal arc at a time, or equivariantly for some
fixed arithmetic Fuchsian group; third and finally, the two essential methods of playing the instrument, one that fixes the
tuning once and for all and another that allows real-time retuning, whether equivariant or not, while playing.

\subsection{The untuned instrument}\label{plastic1}

We shall consider a rendering of the Farey tesselation $\tau_*$ of $\mathbb D$ on a large touch-sensitive computer
screen, where for the moment we can either {\it tap} an edge $e\in\tau_*$ with a finger to produce a tone or {\it hold} it at some point
$p\in e$ to sustain a tone.

Each $e\in\tau_*-\{ e_*\}$ comes equipped with its canonical orientation pointing from lower to higher
Farey generation of its endpoints, extended to $e_*$ pointing from ${0\over 1}$ to ${1\over 0}$ by convention.
Each such edge $e\in\tau_*$ also comes equipped with
its distinguished point $p_0\in e$ lying in two Farey horocycles, that is, equidistant to the Farey horocycles centered at
the endpoints of $e$. 
This point in turn determines the net of points $p_j\in e$, $j \in {\mathbb Z}$, called {\it frets}, at hyperbolic distance ${j\over 2}{{\rm exp}({1\over 2})}$ from the distinguished fret $p_0$ along $e$.

Fix  any edge $e\in\tau_*$. In this untuned setting, tapping $e$ always produces a burst at
27.5 Hz by convention, which is likewise produced continuously when the distinguished fret $p_0\in e$ is held.  
The fret $p_i$ can be tapped or held and emits the note which is $i$ hemitones from that of $p_0$ in whatever tempering seems 
salutary, like keys on a piano keyboard but also like frets on a guitar.

More generally, holding an arbitrary point on $e$ determines a continuum of
frequency response that interpolates that of the frets in any desired manner, so tapping or holding nearby a fret modulates the tempering.  It is this more general plasticity of tone, akin to
a violin or mouth harmonica, that is reflected in the name of the new instrument 

Let us indicate orientation with a color spectrum on each edge, with tick marks for the frets, and a special
$\times$-mark for the distinguished fret.  As with the original hormonica, the timbre of each sound (its character
as determined by its spectrum) and its attack (how the tone itself initiates and decays) can be 
set however one desires, for instance like a staccato saxophone
in one region of $\tau_*$, like a slide guitar in another region, and like a snare drum in yet another.  In principle but impracticably,
any piece of music could be played using one edge of $\tau_*$ for each instrument. 

As in the original hormonica, let us also allow tapping and holding
a complementary ideal triangle to actuate the three tones associated with the distinguished frets in its
three frontier edges, which agree for the untuned $\tau_*$.

\subsection{Tuning the instrument} \label{plastic2} Start from the untuned $\tau=\tau_*$, and recursively suppose a tesselation $\tau$ of $\mathbb D$ is displayed whose vertices are given by the Farey rationals and whose colored edges come equipped with a net of frets at distance ${1\over 2}{{\rm exp}({1\over 2})}$ from one another, together with a distinguished fret denoted $\times$ as before.  There is moreover the assignment of a burst or sustained tone to the respective tapping or holding of a point on any edge in $\tau$, as well as triples of tones of the distinguished frets in the frontier of a complementary ideal triangle for tapping or holding  any point in an ideal triangle complementary to $\tau$.

Given an edge $e\in\tau$, imagine a new type of tap called a {\it pedal-tap}, for instance activated by holding down a foot-pedal
while tapping.  The effect of a pedal-tap on $e\in\tau$ is first of all to perform a flip upon this edge of $\tau$ so as to display
the new edge $f$ in the resulting tesselation $\tau'$.  This edge $f\in\tau'$ also runs between Farey rationals and meets the Farey decoration at points whose hyperbolic distance $\delta$ satisfies ${\rm exp}\,\delta/2\in{\mathbb Z}$ according to Corollary
\ref{cor:int}, and thus extends to a system of frets on $f$ between the two horocycles, where the distinguished fret is 
the equidistant point between them.  (Because of the coefficient one half in our specification of distance between frets, the midpoint is always a fret, whether the lambda length is even or odd.)

The frequency associated to tapping $f$ or holding its distinguished fret is given by $27.5\, \xi^\lambda$, where $\lambda$ is the lambda length of $f$ for the Farey decoration, and this tone is sounded as the second effect of a pedal-tap.  A pedal-tap not only retunes $\tau$ to $\tau'$, but also returns a tone based on the resulting new lambda length, so a series of pedal-taps provides not only a series of tuned instruments, but also a series of tones.

Another method of tuning $\tau$ depends upon the specification of a fixed torsion-free $\Gamma<{\rm PSL}_2$
of finite index preserving $\tau$.  In {\it $\Gamma$-equivariant tuning}, the pedal-tap on an edge $e\in\tau$ flips each edge in the
orbit $\Gamma(e)$ to produce another tesselation, which is again preserved by $\Gamma$.  Each edge in $\Gamma(e)$ has the same lambda length, and after a $\Gamma$-equivariant pedal-tap, the instrument sounds a tone of the resulting corresponding common frequency as before.

This completes the recursive definition of the plastic hormonica $\tau$ and its tuning.

\subsection{Playing the instrument}\label{plastic3}

Having tuned the instrument (perhaps equivariantly), the plastic hormonica $\tau$ is played by tapping or pedal-tapping (perhaps equivariantly) edges of $\tau$, by tapping or holding a point of a complementary ideal triangle, or by holding any point of any edge in $\tau$.
A sequence of such actions yields in particular a sequence of tones, and this is how one plays the instrument.  
There are several embodiments of the plastic hormonica depending upon tuning methods
and whether holding is permitted or just tapping.

\bigskip

\noindent{\bf Minimal Simple Embodiment}: The simplest version of the plastic hormonica is to perform a finite sequence of flips once and for all and then simply allow tapping of edges and triangles.  This is essentially the original hormonica of \cite{Pbook}. 
The polygon $P$ in the proof of Lemma \ref{piano} provides a consecutive sequence of hemitones, which could be colored black and white to mimic the piano keyboard.

There is already here as well as below a method of general utility for the automatic play of a tuned plastic hormonica $\tau$.
Namely, a given edge $e\in\tau$ has among its frets of course the distinguished one, but also the two frets which
lie in the horocycles decorating the endpoints of $e$.  Given a point $q\in e\cap h$ in such a horocycle $h$, we can traverse $h$
with unit hyperbolic distance from $q$ in a suitable unit of time, thereby meeting a sequence $q_k=h\cap e_k$ of points in $h\cap\tau$, for edges $e_k\in\tau$.  The sequence of tones for the lambda lengths of $e_k$, played in this time signature determined by hyperbolic distance along the horocycle, provides a kind of arpeggio automatic-play of the instrument from a specified point of the Farey decoration lying in an edge of $\tau$.

\bigskip

\noindent {\bf Minimal Equivariant Embodiment}: Take any punctured arithmetic surface $F={\mathbb D}/\Gamma$,
so $\Gamma<{\rm PSL}_2$ preserves the Farey tesselation $\tau_*$.  Only sequences of $\Gamma$-equivariant
tunings, that is, finite sequences of pedal-taps, are permitted in this embodiment.  As explained in \cite{Pdec,Pbook}, such sequences correspond to paths in the decorated Teichm\"uller space of $F$.  (These integral lambda lengths correspond to the ``centers of top-dimensional cells,'' flips
to ``crossing codimension-one faces between these cells''
and the transitivity in Remark \ref{jhc} to general position in the path connected decorated Teichm\"uller space.)   

Moreover, the mapping class $\varphi:F\to F$ corresponds to the
sequence of flips from an ideal triangulation $\tau$ of $F$ to $\varphi^{-1}(\tau)$, hence to a sequence of $\Gamma$-equivariant
retunings manifest as a periodic sequence of pedal-taps with its corresponding melody.   

 
There are endless attendant questions and tasks including:  

\leftskip .5in\rightskip .3in

\noindent $\bullet$ Play the arpeggio of each puncture in some examples of equivariant tuning.  Play some periodic homeomorphisms of punctured surfaces.  How does a pseudo-Anosov homeomorphism sound?
How do these several sonic tasks behave under finite-sheeted  covers, possibly branched over the (missing)
punctures as with the punctured solenoid? Which integrally weighted pentagons occur universally, or equivariantly, and how do their pentagon relations sound?

\noindent $\bullet$ Choose any fixed piece of music and ask for the minimal genus ${\mathbb D}/\Gamma$, for $\Gamma<{\rm PSL}_2$ so that the piece could be played on a $\Gamma$-equivariant tuning.  The hyperfan piano Example \ref{piano} shows that any single-voiced tune           
has finite minimal genus bounded just in terms of the number of octaves it spans.  Chords are another matter entirely, as has been discussed.   

\noindent $\bullet$ Is there a natural tempering of the plastic hormonica more sonorous than equal temperament, where the frequency of an edge perhaps depends not only multiplicatively on its integral lambda length as with the ordinary hormonica, but also on the rational Farey labeling of its endpoints?  Could the Minkowski question mark function
(which maps the dyadic tesselation to the Farey tesselation, cf. \cite{Pbook}) play a role in sonorous tempering?

\leftskip=0ex\rightskip=0ex

\bigskip


The two minimal embodiments naturally combine: after a finite sequence of possibly $\Gamma$-equivariant retunings from $\tau_*$, play the resulting instrument with tapping in one embodiment and with $\Gamma$-equivariant pedal-tapping as well in another, providing dynamic real-time perhaps equivariant instrument retuning.

\bigskip 

\noindent {\bf Most Flexible Embodiment}: On top of the full tapping and possibly equivariant
pedal-tapping instrument, we finally include holding points on the edges, to unleash a continuum
of sound for each edge.  Tuples of points can be simultaneously held to produce any musical chord.  This embodiment involves full functionality of the frets,
now with arpeggio auto-play also enabled for each oriented edge.  This is the full instrument which seems to
have the desired nuanced responsiveness for musical expression.

\section{Closing remarks}\label{outtahere}

This overall octave scaling $N=4$ to $27.5\, \xi^{\lambda}=440\, \xi^{\lambda-12N}$ of frequency 
from lambda length was unrealistic and metaphorical,
where $N$ or more generally some other lowest reasonably audible frequency
can be adjusted as part of pre-tuning.  The distinguished fret illustrates
where some fixed octave of the note for the string is located, which is useful for holding edges.

In probing moduli spaces and mapping class groups with the minimal equivariant embodiment, any assignment of
frequency to lambda length will suffice, for instance that of the original hormonica.   
The frequency tuning of the plastic hormonica is specified
simply in order to mimic familiar instruments for musical expression.

A natural extension of the plastic hormonica, among others, allows tapping or holding points of a complementary triangle to activate
not just the single triple of frequencies associated to its frontier edges, 
but different modulations and combinations of this triple depending on the point of contact within
the triangle.

One can imagine circumstances where the discrepancies between our own at least locally Euclidean world and the hyperbolic plane may make
instrument play impractical at the scale of fingers on pictures on screens.  The solution is a renormalization by ${\rm PSL}_2={\rm PSL}_2({\mathbb Z})$
or ${\rm PSL}_2({\mathbb R})$ if the scale becomes problematic, perhaps with several simultaneous screenshots of different regions at different scales, or better yet,
a multi-scale virtual reality implementation, as long as we are being speculative.

One might also dream of an implementation of the minimal equivariant embodiment on a physical model of a punctured surface,
i.e., a higher-genus touch-sensitive screen as the boundary of a body in space, a physical three-dimensional version
of the equivariant plastic hormonica of some fixed topological type.  Again in principle, this could presumably be executed in virtual reality.
Another more pedestrian solution to the difficulties mentioned in the previous paragraph in this equivariant setting would be to display only
the horocycles in the quotient surface, together with their intersection points with the edges of the tesselation, which can be tapped or pedal-tapped.

At the other extreme, there is a purely physical way to imagine the plastic hormonica with fixed tuning as a collection of strings in the hyperbolic
plane with endpoints on the horocycles of the Farey decoration.  What is the actual physics of the sounds produced
by plucking these strings in hyperbolic space?  Does sound perception differ in spaces of different curvatures?

Plans are afoot to implement one or another of these embodiments, or at least listen to some mapping classes.  It will surely be fun, and we can always retreat to the {\sl Mathematician's Last Refuge}: this all could be a little bit too much, I admit, and if so, I hope you agree that at least there seem to be  some interesting questions.

\end{document}